\input amstex\documentstyle {amsppt}  
\pagewidth{12.5 cm}\pageheight{19 cm}\magnification\magstep1
\topmatter
\title Homomorphisms of the alternating group $\Cal A_5$ into reductive groups
\endtitle
\dedicatory{Dedicated to Robert Steinberg on his 80-th birthday}
\enddedicatory
\author G. Lusztig\endauthor
\address Department of Mathematics, M.I.T., Cambridge, MA 02139\endaddress
\thanks Supported in part by the National Science Foundation\endthanks
\endtopmatter
\document

\define\fff{{}_{235}}
\define\bF{\bold F}

\redefine\T{\times}
\redefine\t{\tau}

\define\n{\notin}

\define\a{\alpha}

\define\e{\emptyset}
\redefine\i{^{-1}}
\define\g{\gamma}

\define\m{\mapsto}
\define\x{\boxed}

\define\sh{\sharp}

\define\do{\dots}

\define\sqc{\sqcup}

\define\lra{\leftrightarrow}

\define\sub{\subset}

\define\ti{\tilde}
\define\nl{\newline}
\define\fra{\frac}

\define\ot{\otimes}
\define\bbq{\bar{\bq}_l}

\define\Ad{\text{\rm Ad}}
\define\Hom{\text{\rm Hom}}

\define\sgn{\text{\rm sgn}}
\define\tr{\text{\rm tr}}

\define\he{\heartsuit}

\define\spa{\spadesuit}
\define\di{\diamond}
\define\opl{\oplus}

\define\de{\delta}
\define\ep{\epsilon}
\define\io{\iota}

\define\rh{\rho}
\define\si{\sigma}
\define\th{\theta}
\define\ka{\kappa}
\define\la{\lambda}
\define\ph{\phi}
\define\ps{\psi}

\define\vt{\vartheta}

\define\De{\Delta}

\define\boc{\bold c}

\define\bc{\bold C}

\define\bq{\bold Q}

\define\bv{\bold V}

\define\bz{\bold Z}
\define\bx{\bold X}
\define\ca{\Cal A}

\define\cc{\Cal C}

\define\ce{\Cal E}
\define\cf{\Cal F}
\define\cg{\Cal G}

\define\co{\Cal O}

\define\car{\Cal R}

\define\ct{\Cal T}

\define\cx{\Cal X}
\define\cy{\Cal Y}

\define\fg{\frak g}

\define\bZ{\bar Z}

\define\tG{\ti G}
\define\tH{\ti H}
\define\AL{A}
\define\BL{BL}
\define\DL{DL}
\define\FR{F1}
\define\FRR{F2}
\define\FRRR{F3}
\define\GE{C}
\define\LU{L}
\define\SE{S}

\head Introduction \endhead
\subhead 0.1\endsubhead
Let $\ca_5$ be the group defined by generators $x_2,x_3,x_5$ and relations
$x_2^2=x_3^3=x_5^5=x_2x_3x_5=1$. It is well known that $\ca_5$ is isomorphic to
the alternating group in $5$ letters. Let $k$ be an algebraically closed field
of characteristic $p$ where $p=0$ or a prime $\ge 7$. Let $G$ be a connected
reductive algebraic group over $k$. We shall be interested in the 
classification of homomorphisms $\ca_5@>>>G$ up to conjugation by $G$. When
$G=GL_n$ this is the same as the classification of $n$-dimensional 
representations of $\ca_5$ over $k$ up to isomorphism; this can be reduced to
the classification of irreducible representations of $\ca_5$ which is classical
(there are five of them up to isomorphism, of degrees $1,3,3,4,5$). In the 
general case, there is again something analogous to the notion of irreducible 
representation, which we call a {\it regular} homomorphism $\ca_5@>>>G$. This 
is, by definition, a homomorphism $\ca_5@>>>G$ whose image is not contained in
a Levi subgroup of a proper parabolic subgroup of $G$. Again the classification
of homomorphisms $\ca_5@>>>G$ up to conjugacy can be more or less reduced to 
the analogous question for regular homomorphisms (for $G$ and also for smaller
groups). 

D.D.Frey \cite{\FR},\cite{\FRR},\cite{\FRRR} has classified up to
conjugacy the homomorphisms (resp. non-regular homomorphisms) $\ca_5@>>>G$ for
$G$ of type $E_6,E_7$ (resp. $E_8$) over $\bc$.
For $G$ of type $E_8$ he showed that there
is at least one regular homomorphism $\ca_5@>>>G$, but the problem of 
classifying up to conjugacy the regular homomorphisms $\ca_5@>>>G$ remained 
open. 

Serre \cite{\SE} has suggested that this problem could be attacked using the 
complex representation theory of $E_8(\bF_q)$, $\bF_q$ a finite field. He 
showed that the number $d$ of conjugacy classes of regular homomorphisms 
$\ca_5@>>>E_8(k)$ can be extracted from a certain sum over all irreducible 
characters of $E_8(\bF_q)$, if this sum can be computed. 

The purpose of this paper is to show that this sum can be computed with enough
precision so that a solution of the above problem is obtained. (See Theorem 
4.5.) In fact, we show that $d=1$. Thus, our results, in conjunction with the
earlier work of Frey, complete the classification up to $G$-conjugacy of 
homomorphisms $\ca_5@>>>G$ for $G$ adjoint simple.

One of the key ingredients in our proof is a
collection of inequalities (see 2.2) involving representations of Weyl groups.
While these inequalities are conjectural in general, enough of them can be 
verified (with the aid of a computer) so that the proof goes through. 

I want to thank R.Griess for introducing me to this problem. I also want to
thank J.-P.Serre for making \cite{\SE} available to me and for some useful
discussions.

\head 1. Examples\endhead
\subhead 1.1\endsubhead
If $G$ is a group we denote by $Z_G$ the centre of $G$ and by $G_{der}$ the 
derived group of $G$; we set $G_{ad}=G/Z_G$. If $g\in G$ let $Z_G(g)$ be the 
centralizer of $g$ in $G$ and let $\cc_G(g)$ be the conjugacy class of $g$ in 
$G$. If $H$ is a subgroup of $G$ let $Z_G(H)$ (resp. $N_G(H)$) be the 
centralizer (resp. normalizer) of $H$ in $G$. If $G,G'$ are groups, let
$\Hom(G,G')$ be the set of group homomorphisms $G@>>>G'$.

If $G$ is an affine algebraic group over $k$ (as in 0.1), let $G^0$ be the 
identity component of $G$ and let $\bar G=G/G^0$. 

For a finite set $Y$ let $\x{Y}$ or $\sh(Y)$ be the cardinal of $Y$.

\subhead 1.2\endsubhead
Assume now that $G$ is connected, reductive algebraic group over $k$ of simply
laced type. Let 
$\pi:\tG@>>>G$ be the simply connected covering of $G_{der}$. Let $\ct(G)$ be
the set of maximal tori in $G$. Let $r(G)$ be the rank of $G$ and let $\nu(G)$
be half the number of roots of $G$.

Since $\ca_5$ is perfect, any homomorphism $\ps:\ca_5@>>>G$ has image contained
in $G_{der}$. This gives us a bijection between the set of $G$-orbits on 
$\Hom(\ca_5,G)$ and the set of $G_{der}$-orbits on $\Hom(\ca_5,G_{der})$.
Now, the obvious map from $G$-orbits on $\Hom(\ca_5,G)$ to $G_{ad}$-orbits on
$\Hom(\ca_5,G_{ad})$ is injective. 

(It is enough to show that, if $\ps,\ps'\in\Hom(\ca_5@>>>G)$ are such that 
$h\ps=h\ps'$ where $h:G@>>>G_{ad}$ is the obvious map, then $\ps=\ps'$. For
$x\in\ca_5$ we have $\ps'(x)=\ps(x)z(x)$ where $z\in\Hom(\ca_5,Z_G)$. Since
$\ca_5$ is perfect and $Z_G$ is abelian, we have $z(x)=1$ for all $x$, hence 
$\ps=\ps'$.)

\proclaim{Lemma 1.3 (Serre \cite{\SE})} Let $\ps\in\Hom(\ca_5,G)$. We have 

$\sum_{n=2,3,5}\dim Z_G(\ps(x_n))=\dim(G)+2\dim Z_G(\ps(\ca_5))$.
\endproclaim
This is equivalent to its Lie algebra analogue
$$\sum_{n=2,3,5}\dim\fg^{\ps(x_n)}=\dim\fg+2\dim\fg^{\ca_5}$$
where $\fg=\text{\rm Lie }G$, and the upper index denotes the fixed point set 
of the appropriate automorphism or group of automorphisms of $\fg$. More 
generally, this holds when $\fg$ is replaced by a finite dimensional $k$-vector
space $V$ with an $\ca_5$-action; we may assume that $V$ is an irreducible 
representation of $\ca_5$ (we have complete reducibility since 
$\fra{1}{30}\in k$) and in that case the equality is checked by direct 
computation.

\proclaim{Lemma 1.4} If $\ps\in\Hom(\ca_5,G)$, then $Z_G(\ps(\ca_5))$ is a 
reductive group.
\endproclaim
Let $\cc(n)=\cc_G(\ps(x_n)),n=2,3,5$. Then $\cc(n)$
is an affine variety since $\ps(x_n)$ is a
semisimple element. Consider the affine variety 
$$\bx=\{(g_2,g_3,g_5)\in\cc(2)\T\cc(3)\T\cc(5);g_2g_3g_5=1\}.$$
The group $G$ acts on $\bx$ by conjugation on all factors. By 1.3, the 
dimension of the $G$-orbit of $(g_2,g_3,g_5)\in\bx$ equals 
$$\dim(G)-(\sum_{n=2,3,5}\dim Z_G(g_n)-\dim(G))/2=
\sum_{n=2,3,5}\dim\cc(n)/2;$$
in particular it is independent of the choice $(g_2,g_3,g_5)$. Since any 
$G$-orbit of minimum dimension must be closed in $\bx$, it follows that any 
$G$-orbit in $\bx$ is closed in $\bx$; hence it is affine. By a known criterion
(Richardson) it follows that the isotropy group in $G$ of any point of $\bx$ is
reductive. The lemma is proved.

\proclaim{Lemma 1.5} The following four conditions for $\ps\in\Hom(\ca_5,G)$ 
are equivalent:

(i) $\dim Z_G(\ps(\ca_5))=\dim Z_G$;

(ii) $Z_G(\ps(\ca_5))/Z_G$ is finite;

(iii) any subtorus of $Z_G(\ps(\ca_5))$ is contained in $Z_G^0$;

(iv) $\nu(G)-\sum_{n=2,3,5}\nu(Z_G^0(\ps(x_n)))=r(G_{ad})$.
\endproclaim
This follows from Lemmas 1.3, 1.4.

\subhead 1.6\endsubhead
We say that $\ps\in\Hom(\ca_5,G)$ is {\it regular} if the conditions of Lemma
1.5 are satisfied. Let 
$\Hom_{reg}(\ca_5,G)=\{\ps\in\Hom(\ca_5,G);\ps\text{ regular}\}$.
Note that $\ps\in\Hom(\ca_5,G)$ is regular if and only if its image in 
$\Hom(\ca_5,G_{ad})$ is regular.

We show that the classification of homomorphisms $\ps:\ca_5@>>>G$ up to 
$G$-conjugacy can be essentially reduced to the analogous problem for regular 
homomorphisms. Let $\ps\in\Hom(\ca_5,G)$. Let $S$ be a maximal torus of 
$Z_G(\ps(\ca_5))$ and let $L=Z_G(S)$. (A Levi subgroup of a parabolic subgroup
of $G$.) Using 1.5(iii), we see that $\ps$ defines a regular homomorphism
$\ca_5@>>>L$. This gives us a bijection between the set of $G$-orbits on 
$\Hom(\ca_5,G)$ and the disjoint union over all $G$-conjugacy classes of Levi 
subgroups $L$ of the sets of $N_G(L)$-orbits on $\Hom_{reg}(\ca_5,L)$.

\subhead 1.7\endsubhead
A $\fff$-triple for $G$ is a triple $(C_2,C_3,C_5)$ of conjugacy classes in 
$G$ such that, if $g_n\in C_n$, then $g_n^n=1$, $n=2,3,5$. The {\it type} of 
a 
$\fff$-triple is by definition the sequence $X_2,X_3,X_5$ where $X_n$ is the 
type of the root system of $Z_G^0(g_n)$. Clearly, $G$ has only finitely many 
$\fff$-triples. As pointed out in \cite{\SE}, of particular interest are the 
$\fff$-triples $(C_2,C_3,C_5)$ in $G$ such that
$$\nu(G)-\sum_{n=2,3,5}\nu(Z_G^0(g_n))=r(G_{ad})$$ 
where $g_n\in C_n$. Such $\fff$-triples are called {\it regular}. If 
$\ps\in\Hom_{reg}(\ca_5,G)$ and $C_n=\cc_G(\ps(x_n))$, $n=2,3,5$ then 
$(C_2,C_3,C_5)$ is a regular $\fff$-triple. 

Let $\io$ be an automorphism of order $2$ of $\ca_5$ such that $\io(x_2)$ is 
conjugate to $x_2$, $\io(x_3)$ is conjugate to $x_3$, $\io(x_5)$ is conjugate 
to $x_5^2$. For any $\ps\in\Hom(\ca_5,G)$ let ${}^\io\ps=\ps\circ\io$. If 
$(C_2,C_3,C_5)$ is a regular $\fff$-triple of $G$, then
${}^\io(C_2,C_3,C_5)=(C_2,C_3,C'_5)$, where $C'_5=\{g^2;g\in C_5\}$, is again a
regular $\fff$-triple of $G$.

\subhead 1.8\endsubhead
Assume that $G=G_{ad}$ is simple of type $A_{m-1}$. If $m\in\{2,3\}$ then $G$ 
has exactly two regular $\fff$-triples (interchanged by $\io$) and there are 
exactly two $\ps\in\Hom_{reg}(\ca_5,G)$ up to conjugacy (interchanged by 
$\io$). If $m=4$, then $G$ has exactly one regular $\fff$-triple and there are
exactly two $\ps\in\Hom_{reg}(\ca_5,G)$ up to conjugacy. If $m\in\{5,6\}$ then
$G$ has exactly one regular $\fff$-triple and there is exactly one 
$\ps\in\Hom_{reg}(\ca_5,G)$ up to conjugacy. If $m\ge 7$, then $G$ has no 
regular $\fff$-triple and therefore $\Hom_{reg}(\ca_5,G)=\e$.

The type of a regular $\fff$-triple is $(\e,\e,\e)$ (if $m=2$),$(A_1,\e,\e)$ 
(if $m=3$), $(A_1^2,A_1,\e)$ (if $m=4$), $(A_2A_1,A_1^2,\e)$ (if $m=5$), 
$(A_2^2,A_1^3,A_1)$ (if $m=6$).

For $m=3,5$ (resp. $m=2,6$), $\ps$ comes from an irreducible $m$-dimensional
representation of $\ca_5$ (resp. of the duble cover $SL_2(\bF_5)$ of $\ca_5$, 
regarded as a projective representation of $\ca_5$). For $m=4$, one $\ps$ 
comes from an irreducible $4$-dimensional representation of $\ca_5$, the other
$\ps$ comes from an irreducible $4$-dimensional representation of 
$SL_2(\bF_5)$ which does not factor through $\ca_5$.

\subhead 1.9\endsubhead
Assume that $G=G_{ad}$ is simple of type $D_m$. 

(a) If $m=4$, then $G$ has exactly four regular $\fff$-triples: two of type 
\linebreak $(A_1^4,A_1^3,\e)$ (interchanged by $\io$) and two of type 
$(A_1^4,A_2,\e)$ (interchanged by $\io$); correspondingly, it has exactly four
$\ps\in\Hom_{reg}(\ca_5,G)$ up to conjugacy: two of them (interchanged by 
$\io$) come from $8$-dimensional orthogonal representations of $\ca_5$ which 
decompose as $1+3+4$ into irreducibles and two of them (interchanged by $\io$)
come from $8$-dimensional orthogonal representations of $\ca_5$ which decompose
as $3+5$ into irreducibles. 

(b) If $m=5$, then $G$ has exactly one regular $\fff$-triple; it has type
\linebreak $(A_3A_1^2,A_2A_1^2,A_1^2)$; also $G$ has exactly two 
$\ps\in\Hom_{reg}(\ca_5,G)$ up to conjugacy; one of them comes from a 
$10$-dimensional orthogonal representation of $\ca_5$ which decomposes as 
$3+3+4$ into irreducibles; the other comes from a $10$-dimensional orthogonal 
representation of $\ca_5$ which decomposes as $1+4+5$ into irreducibles. 

(c) If $m=6$, then $G$ has exactly three regular $\fff$-triples; one has type
\linebreak $(A_3^2,A_3A_1^2,A_1^4)$ and the other two (interchanged by $\io$) 
have type \linebreak $(A_3^2,A_3A_1^2,A_2A_1)$; correspondingly, it has exactly
three $\ps\in\Hom_{reg}(\ca_5,G)$ up to conjugacy: one of them comes from a 
$12$-dimensional orthogonal representation of $\ca_5$ which decomposes as 
$1+3+3+5$ into irreducibles; the other two (interchanged by $\io$) come from 
$12$-dimensional orthogonal representations of $\ca_5$ which decompose as 
$3+4+5$ into irreducibles.

(d) If $m=8$, then $G$ has exactly one regular $\fff$-triple; it has type
\linebreak $(D_4A_3,A_4A_3,A_2^2A_1^2)$; it also has exactly one 
$\ps\in\Hom_{reg}(\ca_5,G)$ up to conjugacy: it comes from a $16$-dimensional
orthogonal representation of $\ca_5$ which decomposes as $1+3+3+4+5$ into 
irreducibles.

If $m=7$ or $m>9$ then $G$ has no regular $\fff$-triples; therefore,\linebreak
$\Hom_{reg}(\ca_5,G)=\e$.

We see that in these cases each regular $\ps$ comes from a homomorphism
$\ca_5@>>>SO_{2m}$.

\subhead 1.10\endsubhead
Assume that $G=G_{ad}$ is simple of type $E_m$. From \cite{\FR},
\cite{\FRR},\cite{\FRRR}  we see that:

if $m=6$, then $G$ has exactly two regular $\fff$-triples (interchanged by 
$\io$); they have type $(A_5A_1,A_2^3,A_2A_1^2)$;

if $m=7$, then $G$ has exactly two regular $\fff$-triples (interchanged by 
$\io$); they have type $(A_7,A_5A_2,A_3A_2A_1)$;

if $m=8$, then $G$ has exactly one regular $\fff$-triple; it has type
$(D_8,A_8,A_4^2)$.
\nl
From \cite{\FR},\cite{\FRR},\cite{\FRRR} we see also that in each of these 
cases, a regular $\fff$-triple is associated with at least one 
$\ps\in\Hom_{reg}(\ca_5,G)$. (Another proof of this fact is given by Lemma 
4.3.) This $\ps$ is in fact unique up to conjugacy. (See \cite{\FRR},
\cite{\FRRR} for $E_6,E_7$ and $p=0$; see 4.5, 4.6 for the general case.)

\head 2. Inequalities\endhead
\subhead 2.1\endsubhead
Let $W$ be a Weyl group and let $V$ be the reflection representation of $W$ 
(over $\bc$). Let $\bv=\opl_{j\ge 0}\bv_j$ be the algebra of polynomials 
$V@>>>\bc$ modulo the ideal generated by the $W$-invariant polynomials of 
degree $>0$; here $\bv_j$ is the image of the space of polynomials of degree
$j$. Now $\dim\bv=\x{W}$ and $W$ acts naturally on $\bv$ preserving the 
grading. Let $\ce(W)$ be the set of irreducible representations of $W$ over 
$\bc$ (up to isomorphism). For any $E\in\ce(W)$ let $b_W(E)$ (resp. $b'_W(E)$)
be the smallest (resp. largest) integer $j$ such that $\Hom_W(E,\bv_j)\ne 0$.
Note that $b_W(E),b'_W(E)$ are well defined.

Assume now that $W$ is the Weyl group of $G$ as in 1.2. Let $\nu=\nu(G)$. Let 
$(C_2,C_3,C_5)$ be a $\fff$-triple for $G$. Let 
$(g_2,g_3,g_5)\in C_2\T C_3\T C_5$. We denote by $W_n$ (resp. by $\nu_n$) the 
Weyl group (resp. the number of positive roots) of $Z_G^0(g_n)$, $n=2,3,5$. We
may regard $W_2,W_3,W_5$ as Weyl subgroups of $W$. Let 
$\bv^{(n)}=\opl_j\bv^{(n)}$ be defined in terms of $W_n$ in the same way as
$\bv=\opl_j\bv_j$ is defined in terms of $W$. For any $E\in\ce(W)$ let $b_n(E)$
(resp. $b'_n(E)$) be the smallest (resp. largest) integer $j$ such that 
$\Hom_{W_n}(E|_{W_n},\bv_j^{(n)})\ne 0$. Note that $b_n(E),b'_n(E)$ are well 
defined for $n=2,3,5$.

Let $\sgn$ be the sign representation of $W$. We have
$$b(E)=\nu-b'(E\ot\sgn),\quad b_n(E)=\nu_n-b'_n(E\ot\sgn),\quad
 n=2,3,5.\tag a$$
Now $\ce(W)$ is a union of equivalence classes called {\it families}, see
\cite{\LU, 4.2}. For example the unit representation $1$ is a family by itself.
If $\cf$ is a family, then $\{E\ot\sgn; E\in\cf\}$ is again a family. For any 
family $\cf$ of $W$ we set
$$a(\cf)=\min_{E\in\cf}b(E),\quad a'(\cf)=\max_{E\in\cf}b'(E),$$
$$a_n(\cf)=\min_{E\in\cf}b_n(E),\quad a'_n(\cf)=\max_{E\in\cf}b'_n(E).$$
From (a) we deduce that
$$a(\cf)=\nu-a'(\cf\ot\sgn),\quad a_n(\cf)=\nu_n-a'_n(\cf\ot\sgn),\quad
n=2,3,5.\tag b$$

\proclaim{Proposition 2.2} In the setup of 2.1, assume that either 

(1) $(C_2,C_3,C_5)$ is regular, or 

(2) any simple component of $G_{ad}$ has rank $\le 8$ and is not of type $E_8$.
\nl
Then for any family $\cf$ we have

(a) $a(\cf)-\sum_{n=2,3,5}a_n(\cf)\le r(G_{ad})$;
\nl
equivalently (see 2.1(b)), 

(b) $\sum_{n=2,3,5}a'_n(\cf)-a'(\cf)\le-\nu+\sum_n\nu_n+r(G_{ad})$.
\nl
Moreover, if $(C_2,C_3,C_5)$ is not regular, then (a) and (b) are strict
inequalities. If $(C_2,C_3,C_5)$ is regular, then (a) is strict for
$\cf\ne\{\sgn\}$ and (b) is strict for $\cf\ne\{1\}$. If $(C_2,C_3,C_5)$ is 
regular, then (a) is an equality for $\cf=\{\sgn\}$ and (b) is an equality for
$\cf=\{1\}$. 
\endproclaim
This can be checked with the aid of a computer. We may assume that $G=G_{ad}$ 
is simple. In the case where $G$ is of type $E_8$ and $(C_2,C_3,C_5)$ is
regular, I have used tables in \cite{\BL}, \cite{\AL}. If we are not in this
case then $G$ has rank $\le 8$ and is not of type $E_8$. Then, instead of  
\cite{\BL}, \cite{\AL}, one can use the CHEVIE package \cite{\GE} to get tables
for the quantities $b_n(E)$ (or rather the analogous quantities where $W_n$ is
replaced by any reflection subgroup of $W$). From this information, one can 
check the required inequalities by hand. (I am grateful to M. Geck for his help
with the CHEVIE package.)

\subhead 2.3\endsubhead
We expect that 2.2 continues to hold even if we drop the assumptions (1),(2).
One can check that 2.2(a),(b) hold for $G$ of type $A$ without the assumptions
(1),(2). When $(C_2,C_3,C_5)$ is  regular, then 2.2(a) for $\cf=\{\sgn\}$ is 
just the equality $\nu-\sum_{n=2,3,5}\nu_n=r(G_{ad})$.

\subhead 2.4\endsubhead
We illustrate 2.2(a) in the case when $G=G_{ad}$ is of type $E_8$ and 
$(C_2,C_3,C_5)$ is regular. The subgroups $W_n$ of $W$ are as in 1.10. In the 
following list, there is one line for each family $\cf$ of $W$, containing the
information

$D; a(\cf)-a_2(\cf)-a_3(\cf)-a_5(\cf)=?$
\nl
where $D$ is the degree of the "special" representation in $\cf$. Note that 
$\cf$ is determined by $D$ and $a(\cf)$. 

$1; 0-0-0-0=0$

$8; 1-1-1-1=-2$

$35; 2-2-1-0=-1$

$112; 3-0-0-0=3$

$210; 4-0-1-0=3$

$560; 5-1-1-0=3$

$567;  6-2-1-1=2$

$700; 6-1-0-0=5$

$1400; 7-2-0-0=5$

$1400; 8-2-0-0=6$

$3240; 9-3-1-0=5$

$2268; 10-4-2-0=4$

$2240; 10-3-1-0=6$

$4096; 11-4-2-1=4$

$525; 12-4-3-1=4$

$4200; 12-4-2-0=6$

$2800; 13-4-3-1=6$

$4536; 13-5-2-0=6$

$2835; 14-6-2-0=6$

$6075; 14-6-3-1=4$

$4200; 15-5-3-1=6$

$5600; 20-8-5-2=5$

$4480; 16-6-3-0=7$

$5600; 21-8-5-2=6$

$4200; 21-9-5-2=5$

$6075; 22-10-5-2=5$

$2835; 22-10-5-2=5$

$4536; 23-9-5-2=7$

$2800; 25-11-7-3=4$

$4200; 24-10-6-2=6$

$525;  36-16-9-5=6$

$4096; 26-11-7-3=5$

$2240; 28-12-7-3=6$

$2268; 30-12-8-4=6$

$3240; 31-13-8-4=6$

$1400; 32-14-9-5=4$

$1400; 37-16-9-5=7$

$700; 42-18-12-6=6$

$567; 46-22-13-7=4$

$560; 47-21-13-7=6$

$210; 52-24-16-8=4$

$112; 63-28-18-10=7$

$35; 74-34-22-12=6$

$8; 91-43-28-16=4$

$1; 120-56-36-20=8$.

\head 3. $\bF_q$-structures \endhead
\subhead 3.1\endsubhead
In this section we fix $G$ as in 1.2. We assume that $p$ is a prime $\ge 7$ and
that $k$ is an algebraic closure of the field with $p$ elements. Let $q$ be a 
power of $p$ and let $\bF_q$ be the subfield of $k$ with $q$ elements. Let 
$F:G@>>>G$ be the Frobenius map corresponding to an $\bF_q$-rational structure
on $G$. 

For any integer $m\ge 1$ we write $m=m_pm_{p'}$ where $m_p$ is a power of $p$
and 
$m_{p'}$ is prime to $p$. If $\ph$ is a map of a set $Y$ into itself, we set
$Y^\ph=\{y\in Y;\ph(y)=y\}$.

\subhead 3.2\endsubhead
Let $\ep_G=(-1)^{\bF_q-\text{\rm rank of } G}$. 
For any $T\in\ct(G)^F$ and any 
$\vt\in\Hom(T^F,\bc^*)$ let $R_{T,G}^\th:G^F@>>>\bc$ be the character of the 
virtual representation denoted in \cite{\DL, 4.3} by $R_T^\th$. (We choose a 
prime number $l,l\ne p$, an algebraic closure $\bar Q_l$ of the $l$-adic
numbers and a field isomorphism $\bar Q_l\cong\bc$; then $R_T^\th$, a virtual 
representation over $\bbq$, can be viewed as a virtual representation over
$\bc$.) Let $\ce(G^F)$ be the set of irreducible (complex) characters of $G^F$.
Let $\ce_1(G^F)$ be the set of irreducible unipotent characters of $G^F$ that 
is, the irreducible characters $\rh$ of $G^F$ such that 
$(\rh:R_{T,G}^1)_G\ne 0$ for some $T\in\ct^F$. Here $(:)_G$ is the standard 
hermitian inner product of class functions $G^F@>>>\bc$. Let $\co(G^F)$ be the
set of all class functions $\ph:G^F@>>>\bc$ such that $(\ph:R_{T,G}^\th)_G=0$
for any $T\in\ct(G)^F$ and any $\vt\in\Hom(T^F,\bc^*)$.

\proclaim{Lemma 3.3} Let $\rh\in\ce_1(G^F)$. There exists a unique 
$G^F$-invariant function $f^\rh:\ct(G)^F@>>>\bc$ and a unique $\xi\in\co(G^F)$
such that
$$\rh=\x{G^F}\i\sum_{T\in\ct(G)^F}\x{T^F}f^\rh(T)R_{T,G}^1+\xi.$$
\endproclaim
This follows immediately from the definition of $\ce_1(G^F)$ and from the
orthogonality relations \cite{\DL, 6.8}. 

\proclaim{Lemma 3.4} Let $g\in G^F$ be semisimple. Let $T\in\ct(G)^F, 
\vt\in\Hom(T^F,\bc^*)$. We have
$$R_T^\vt(g)=\ep_{Z^0_G(g)}\ep_T\x{T^F}\i\x{Z_G(g)}_{p'}
\sum_{g'\in\cc_{G^F}(g)\cap T}\vt(g').$$
\endproclaim
See \cite{\DL, 7.2}.

\subhead 3.5\endsubhead
For any $G^F$-invariant function $f:\ct(G)^F@>>>\bc$ let
$$R_{f,G}=\x{G^F}\i\sum_{T\in\ct(G)^F}\x{T^F}f(T)R_{T,G}^1.$$
If $H$ is a connected reductive $F$-stable subgroup of $G$ with $r(H)=r(G)$ and
$f$ is as above, we shall write $R_{f,H}$ instead of $R_{f|_{\ct(H)^F},H}$.

\proclaim{Lemma 3.6} In the setup of 3.5 we have
$\sum_{T\in\ct(G)^F}\ep_G\ep_Tf(T)=\x{G^F}_pR_{f,G}(1)$.
\endproclaim
$R_{f,G}(1)$ may be computed using the formula in \cite{\DL, 7.1} for
$R_{T,G}^1(1)$; the lemma follows.

\subhead 3.7\endsubhead
We denote by $G^F_*$ (resp. $G^F_\di$) the group of all 
$\th\in\Hom(G^F,\bc^*)$ such that $\th|_{\pi(\tG)^F}=1$ (resp. 
$\th|_{\pi(\tG^F)}=1$). Clearly, $G^F_*\sub G^F_\di$.

If $T\in\ct(G)^F$ and $\vt\in\Hom(T^F,\bc^*)$, let $S(T,\vt)$ be the subgroup 
of $G$ defined as in \cite{\DL, 5.19}. If $H$ is an $F$-stable connected 
reductive subgroup of $G$ with $r(H)=r(G)$,
 let $H^F_\spa$ be the set of all $\th\in H^F_\di$ 
such that for some (or equivalently, any) $T\in\ct(H)^F$ we have 
$S(T,\th|_{T^F})=H$. Let $H^F_\he$ be the set of all $\th\in H^F_\di$ such that
$\th=\th'\th''$ for some $\th'\in H^F_\spa,\th''\in H^F_*$.

Let $\cx_{G,F}$ be the set of all pairs $(H,\th)$ where $H$ is an $F$-stable 
connected reductive subgroup of $G$ with $r(H)=r(G)$ and $\th\in H^F_\spa$. 
Let $\cx'_{G,F}$ be
the set of all $F$-stable connected reductive subgroups $H$ of $G$ such that 
$H^F_\spa\ne\e$.

\proclaim{Lemma 3.8} Assume that $Z_G=Z_G^0$. Let $H\in\cx'_{G,F}$. Then any 
orbit of the obvious $N_G(H)^F$-action on $H^F_\spa$ has cardinal 
$\x{N_G(H)^F/H^F}$.
\endproclaim
Let $\th\in H^F_\spa$. Let $S_\th$ be the stabilizer of $\th$ in $N_G(H)^F$. We
show that $S_\th=H^F$. The inclusion $H^F\sub S_\th$ is immediate. Now let 
$g\in S_\th$. We must show that $g\in H^F$. Let $T\in\ct(H)^F$ be maximally
split. Since $gHg\i=H$, we see that $gTg\i\in\ct(H)^F$ is maximally split hence
$gTg\i=hTh\i$ for some $h\in H^F$. Replacing $g$ by $gh\i$ we see that we may 
assume that $gTg\i=T$. Now conjugation by $g$ keeps fixed $\th|_{T^F}$. By
\cite{\DL, 5.13}, $g$ is a product of elements in $N_H(T)$. In particular, 
$g\in H$. Hence $g\in H^F$. The lemma is proved.

\subhead 3.9\endsubhead
Let $(H,\th)\in\cx_{G,F}$. Let $\ce_{H,\th}(G)$ be the set of all 
$\rh\in\ce(G^F)$ such that $(\rh:R_{T,G}^{\th_T})_G\ne 0$ for some 
$T\in\ct(H)^F$. (We write $\th_T$ instead of $\th|_{T^F}$.) According to 
\cite{\DL, 6.3, 5.20}, if $Z_G=Z_G^0$, we have a partition

(a) $\ce(G^F)=\sqc_{(H,\th)}\ce_{H,\th}(G^F)$
\nl
where $(H,\th)$ runs over a set of representatives for the $G^F$-orbits on 
$\cx_{G,F}$. 

According to \cite{\LU}, if $Z_G=Z_G^0$, for any $(H,\th)\in\cx_{G,F}$, there 
exists $\ka\in\{1,-1\}$ and a bijection 

(b) $\ce_1(H^F)\lra\ce_{H,\th}(G^F)$ 
\nl
(denoted by $\rh\lra\rh^{H,\th}$) such that for any $\rh\in\ce_1(H^F)$ and any
$T\in\ct(H)^F$ we have $(\rh:R_{T,H}^1)_H=\ka(\rh^{H,\th}:R_{T,G}^{\th_T})_G$.

\proclaim{Lemma 3.10} Assume that $Z_G=Z_G^0$. Let $(H,\th)\in\cx_{G,F}$. Let 
$\ka$ be as above. Let $\rh\in\ce_1(H^F)$. Let $f^\rh:\ct(H)^F@>>>\bc$ be the
$H^F$-invariant function such that
$\rh=\x{H^F}\i\sum_{T\in\ct(H)^F}\x{T^F}f^\rh(T)R_{T,H}^1+\xi$,
$\xi\in\co(H^F)$ (see 3.3).

(a) We have $\rh^{H,\th}
=\ka\x{H^F}\i\sum_{T\in\ct(H)^F}\x{T^F}f^\rh(T)R_{T,G}^{\th_T}+\xi'$ with
$\xi'\in\co(G^F)$.

(b) Let $g\in G^F$ be semisimple. We have
$$\rh^{H,\th}(g)=\ka\x{Z_G(g)^F}_{p'}\sum_{g'}\ep_{Z^0_G(g)}\ep_{Z^0_H(g')}
\x{Z_H(g')^F}_{p'}\i R_{f^\rh,Z_H^0(g')}(1)\th(g')$$
where $g'$ runs through a set of representatives for the conjugacy classes in
$H^F$ that are contained in $\cc_{G^F}(g)$.

(c) We have $\rh^{H,\th}(1)=\ka\x{G^F}_{p'}\x{H^F}_{p'}\i\ep_G\ep_H
R_{f^\rh,H}(1)$.
\endproclaim
To prove the identity in (a) it is enough to show that both sides have the same
inner product with $R_{T',G}^{\th_{T'}}$ for any $T'\in\ct(H)^F$. By 
\cite{\DL, 6.8}, the inner product of the right hand side with 
$R_{T',G}^{\th_{T'}}$ is
$$\ka\sum_{T\in\ct(H)^F}\fra{f^\rh(T)}{\x{H^F}}\sh
(g_1\in G^F;g_1Tg_1\i=T',\Ad(g_1) \text{ carries $\th_T$ to }\th_{T'})\tag d$$
For any $g_1$ in (d) we see that $\Ad(g_1)$ carries $S(T,\th_T)$ to
$S(T',\th_{T'})$ hence it carries $(H,\th)$ to $(H,\th)$. By the argument in
the proof of 3.8 we have $g_1\in H^F$. Thus, (d) is equal to
$$\ka\x{H^F}\i\sum_{T\in\ct(H)^F}f^\rh(T)\sh(g_1\in H^F;g_1Tg_1\i=T').$$
Now 
$$\align&(\rh^{H,\th}:R_{T',G}^{\th_{T'}})_G=\ka(\rh:R_{T',H}^1)_H\\&=
\ka\x{H^F}\i\sum_{T\in\ct(H)^F}f^\rh(T)\sh(h\in H^F;hTh\i=T').\endalign$$
This proves (a). We prove (b). By \cite{\DL, 7.5} we have $\xi'(g)=0$. Hence
$$\align&\rh^{H,\th}(g)=\ka\x{H^F}\i\sum_{T\in\ct(H)^F}\x{T^F}f^\rh(T)
R_{T,G}^{\th_T}(g)\\&=\ka\x{H^F}\i\sum_{T\in\ct(H)^F}f^\rh(T)
\ep_{Z^0_G(g)}\ep_T\x{Z_G(g)^F}_{p'}\sum_{g'\in\cc_{G^F}(g)\cap T}\th(g')\\&
=\ka\x{H^F}\i\sum_{g'\in\cc_{G^F}(g)\cap H}
\ep_{Z^0_G(g)}\ep_{Z^0_H(g')}\sum_{T\in\ct(H)^F;g'\in T}f^\rh(T)\\&
\T \ep_{Z^0_H(g')}\ep_T\x{Z_G(g)^F}_{p'}\th(g')\\&
=\ka\x{H^F}\i\sum_{g'\in\cc_{G^F}(g)\cap H}\ep_{Z^0_G(g)}\ep_{Z^0_H(g')}
\x{Z^0_H(g')^F}_pR_{f,Z^0_H(g')}(1)\x{Z_G(g)^F}_{p'}\th(g')\endalign$$
and (b) follows. Now (c) is a special case of (b). The lemma is proved.

\medpagebreak

The following result is well known.
\proclaim{Lemma 3.11} Let $A,B,C$ be three conjugacy classes in $G^F$ and let
$(a,b,c)\in A\T B\T C$. We have
$$\sh\{(\ti a,\ti b,\ti c)\in A\T B\T C;\ti a\ti b\ti c=1\}
=\fra{\x{A}\x{B}\x{C}}{\x{G^F}}
\sum_{\rh\in\ce(G^F)}\fra{\rh(a)\rh(b)\rh(c)}{\rh(1)}.\tag a$$
\endproclaim

\subhead 3.12\endsubhead
If $A,B,C$ are three semisimple conjugacy classes in $G^F$, we set
$$\align&\cg_{G,F;A,B,C}=\sum_{\rh\in\ce_1(G^F)}q^{-\nu(Z_G^0(a))-\nu(Z_G^0(b))
-\nu(Z_G^0(c))+\nu(G)-r(G_{ad})}\\&
\T R_{f^\rh,Z_G^0(a)}(1)R_{f^\rh,Z_G^0(b)}(1)R_{f^\rh,Z_G^0(c)}(1)
R_{f^\rh,G}(1)\i\endalign$$
where $f^\rh:\ct(G)^F@>>>\bc$ are defined as in 3.3 and $(a,b,c)\in A\T B\T C$.

If, in addition, $H\in\cx'_{G,F}$ and $A',B',C'$ are three semisimple conjugacy
classes in $H^F$ such that $A'\sub A,B'\sub B,C'\sub C$, we set
$$\align&t_{G,F;H;A,B,C;A',B',C'}=
\ep_{Z^0_G(a')}\ep_{Z^0_H(a')}\ep_{Z^0_G(b')}\ep_{Z^0_H(b')}
\ep_{Z^0_G(c')}\ep_{Z^0_H(c')}\\&\T \fra{\x{A}\x{B}\x{C}}{\x{G^F}}
\fra{\x{Z_G(a)^F}_{p'}\x{Z_G(b)^F}_{p'}\x{Z_G(c)^F}_{p'}}
{\x{Z_H(a')^F}_{p'}\x{Z_H(b')^F}_{p'}\x{Z_H(c')^F}_{p'}}
\fra{\x{H^F}_{p'}}{\x{G^F}_{p'}}\\&
\T q^{-\de+\nu(Z_H^0(a'))+\nu(Z_H^0(b'))+\nu(Z_H^0(c'))-\nu(H)+r(H)}\endalign$$
where $(a',b',c')\in A'\T B'\T C'$ and 
$$\de=3\nu(G)-\nu(Z_G^0(a))-\nu(Z_G^0(b))-\nu(Z_G^0(c)).$$
Note that 
$\ep_{Z^0_G(a)}=\ep_{Z^0_G(a')},\ep_{Z^0_G(b)}=\ep_{Z^0_G(b')},
\ep_{Z^0_G(c)}=\ep_{Z^0_G(c')}$.
We set
$$\align&Y_{G,F;H;A,B,C;A',B',C'}\\&=t_{G,F;H;A,B,C;A',B',C'}\cg_{H,F;A',B',C'}
q^{-\dim Z^0_H}\sum_{\th\in H_\spa^F}\th(a'b'c'),\endalign$$
$$X_{G,F;H;A,B,C}=\sum_{A',B',C'}Y_{G,F;H;A,B,C;A',B',C'},$$
where $A'$ (resp. $B',C'$) runs through the conjugacy classes in $H^F$ that are
contained in $A$ (resp. $B,C$).

Let $\bar\cx'_{G,F}$ be a set of representatives for the orbits of the natural 
$G^F$-action on $\cx'_{G,F}$.

\proclaim{Lemma 3.13} Assume that $Z_G=Z_G^0$. We have
$$q^{-\de}\sh\{(\ti a,\ti b,\ti c)\in A\T B\T C;\ti a\ti b\ti c=1\}=
\sum_{H\in\bar\cx'_{G,F}}\fra{X_{G,F;H;A,B,C}}{\x{N_G(H)^F/H^F}}.$$
\endproclaim
We rewrite the right hand side of 3.11(a) using the partition 
3.9(a), the bijection 3.9(b) and using Lemma 3.8; we obtain
$$\sh\{(\ti a,\ti b,\ti c)\in A\T B\T C;\ti a\ti b\ti c=1\}=
\sum_{H\in\bar\cx'_{G,F}}\x{N_G(H)^F/H^F}\i X'_{G,F;H;A,B,C}$$
where for any $H\in\cx'_{G,F}$ we set
$$X'_{G,F;H;A,B,C}=\fra{\x{A}\x{B}\x{C}}{\x{G^F}}
\sum_{\th\in H_\spa^F}\sum_{\rh\in\ce_1(H^F)}
\fra{\rh^{H,\th}(a)\rh^{H,\th}(b)\rh^{H,\th}(c)}{\rh^{H,\th}(1)}.$$
We evaluate $\rh^{H,\th}(a),\rh^{H,\th}(b),\rh^{H,\th}(c)$ using 3.10(b), and 
$\rh^{H,\th}(1)$ using 3.10(c); we obtain 
$X'_{G,F;H;A,B,C}=q^\de X_{G,F;H;A,B,C}$ (see 3.12). The lemma follows.

\proclaim{Lemma 3.14} (a) In the setup of 3.12 we have 
$|t_{G,F;H;A,B,C;A',B',C'}|\le\boc$ where $\boc>0$ is an integer depending only
on the Coxeter graph on $G$ (and not on $q,F,H,A,B,C,A',B',C'$).

(b) More precisely,
$$\align&|t_{G,F;H;A,B,C;A',B',C'}-\fra{\ep_{Z^0_G(a')}\ep_{Z^0_H(a')}
\ep_{Z^0_G(b')}\ep_{Z^0_H(b')}\ep_{Z^0_G(c')}\ep_{Z^0_H(c')}}
{\x{\bZ_H(a')^F}\x{\bZ_H(b')^F}\x{\bZ_H(c')^F}}|\le\boc'q\i\endalign$$
where $\boc'>0$ is an integer depending only on the Coxeter graph on $G$ (and 
not on $q,F,H,A,B,C,A',B',C'$).
\endproclaim
This is easily checked.

\proclaim{Lemma 3.15} Assume that $Z_G=Z_G^0$.

(a) There exists an integer $\boc>0$ depending only on the Coxeter graph of 
$G$ (and not on $q,F$) such that $\x{\bar\cx'_{G,F}}\le\boc$.

(b) There exists an integer $\boc'>0$ depending only on the Coxeter graph of 
$G$ (and not on $q,F$) such that for any $H\in\cx'_{G,F}$ we have
$\x{H^F_\he-H^F_\spa}\le\boc'q^{\dim(Z_H^0)-1}$.
\endproclaim
Let $G^*$ be a connected reductive group of type dual to that of $G$. Let 
$F':G^*@>>>G^*$ be the Frobenius map corresponding to an $\bF_q$-rational 
structure on $G$. Let $\cy_{G^*,F'}$ be the set of all pairs $(E,g)$ where 
$g\in G^*{}^{F'}$ is semisimple and $E=Z_{G^*}(g)$. Let $\cy'_{G^*,F'}$ be the
set of all subgroups $E$ of $G^*$ such that $(E,g)\in\cy_{G^*,F'}$ for some 
$g$. For $E\in\cy'_{G^*,F'}$ let $Z_{E,\spa}=\{g\in Z_E;Z_{G^*}(g)=E\}$. Let 
$Z_{E,\he}$ be the union of all connected components of $Z_E$ that contain some
point in $Z_{E,\spa}^{F'}$. Let $\bar\cy'_{G^*,F'}$ be a set of representatives
for the orbits of the natural $G^*{}^{F'}$-action on $\cy'_{G^*,F'}$.

Using \cite{\DL, 5.20, 5.24} we can reduce the statements of the lemma to 
statements about $G^*$:

(a${}'$) There exists an integer $\boc>0$ depending only on the Coxeter graph 
of $G^*$ (and not on $q,F'$) such that $\x{\bar\cy'_{G^*,F'}}\le\boc$.

(b${}'$) There exists an integer $\boc'>0$ depending only on the Coxeter graph
of $G^*$ (and not on $q,F'$) such that for any $(E,g)\in\cy'_{G^*,F'}$ we have
$\x{Z_{E,\he}^{F'}-Z^{F'}_{E,\spa}}\le\boc'q^{\dim(Z_E^0)-1}$.
\nl
The proof of (a${}'$) is standard. For (b${}'$) we note that any connected 
component of $Z_{E,\he}$ is an irreducible variety isomorphic to $Z_E^0$ and 
its intersection with $Z_{E,\spa}$ is open and dense in it. The lemma is 
proved.

\subhead 3.16\endsubhead
Let $H$ be an $F$-stable connected reductive subgroup of $G$ such that
$r(H)=r(G)$. Let $W$ be the Weyl group of $G$ and let $W'$ be the Weyl group of
$H$. Note that $F$ acts naturally on $W$ and on $W'$. 

Let $T_0,B_0$ be a pair consisting of an $F$-stable maximal torus of $G$ and an
$F$-stable Borel subgroup of $G$ containing $T_0$. Let $T'_0,B'_0$ be a pair 
consisting of an $F$-stable maximal torus of $H$ and an $F$-stable Borel 
subgroup of $H$ containing $T'_0$. We identify $W=N_G(T_0)/T_0$,
$W'=N_H(T'_0)/T'_0$ in the standard way. We choose $\g\in G$ such that 
$\g T_0\g\i=T'_0$. Let $m=\g\i F(\g)$. We have $m\in N_G(T_0)$. We have an 
imbedding $W'\sub W$ induced by $N_H(T'_0)@>>>N_G(T_0), h\m\g\i h\g$. We
identify $W'$ with its image under this imbedding. The map $F:W'@>>>W'$
corresponds under this identification to the restriction of the map $F':W@>>>W$
induced by $N_G(T_0)@>>>N_G(T_0), g\m mF(g)m\i$.

Let $E$ be an irreducible representation of $W$ with a given isomorphism
$\la_E:E@>>>E$ of finite order such that $\la_E(we)=F(w)\la_E(e)$ for all 
$w\in W,e\in E$. 

Define $f_E:\ct(G)^F@>>>\bc$ by $f_E(T)=\tr(w_T\la_E:E@>>>E)$ where $w_T$ is 
the image in $W$ of $z\i F(z)\in N_G(T_0)$ where $z\in G$ is such that 
$zT_0z\i=T$. (Note that $f_E(T)$ is independent of the choice of $z$.)

Define $\bv'_j$ in terms of $W'$ in the same way as $\bv_j$ was
defined in 2.1 in terms of $W$. Then $W'$ acts naturally on $\bv'_j$ and there
is a natural isomorphism $\la_j:\bv'_j@>>>\bv'_j$ of finite order such that 
$\la_j(w'x)=F'(w')\la_j(x)$ for all $w'\in W',x\in\bv'_j$. 

We define a linear map $\ti\la_j:\Hom_{W'}(E,\bv'_j)@>>>\Hom_{W'}(E,\bv'_j)$ by
$(\ti\la_j(\xi))(e)=\la_j\i\xi(m\la_E(e))$ for all $\xi\in\Hom_{W'}(E,\bv'_j)$,
$e\in E$. (One checks easily that $\ti\la_j(\xi)\in\Hom_{W'}(E,\bv'_j)$.) One 
also checks that $\ti\la_j$ {\it is a map of finite order}. With this notation,
we have
$$R_{f_E,H}(1)=\sum_{j\ge 0}\tr(\ti\la_j,\Hom_{W'}(E,\bv'_j))q^j.\tag a$$
This follows in a standard way from the definitions using the formula 
\cite{\DL, 7.1} for $R_{T,G}^1(1)$.

\proclaim{Lemma 3.17} Let $H\in\cx'_{G,F}$. Let $g\in H^F$ be a semisimple 
element. Assume that $Z^0_H(g)$ is split over $\bF_q$ and that $F$ acts 
trivially on $\bZ_H(g)$. Let $g'\in H^F$ be such that $g,g'$ are conjugate in
$H$. Then:

(a) $\ep_{Z_G^0(g')}=\ep_{Z_H^0(g')}$;

(b) $F$ acts trivially on $\bZ_H(g')$. 
\endproclaim
We can find a maximal torus $T$ of $Z_H^0(g)$ such that $F(T)=T$ and $T$ is 
$\bF_q$-split. Let $T'$ be a maximal torus of $Z_H^0(g')$ such that $FT'=T'$. 
We can find $x\in H$ such that $xgx\i=g',xTx\i=T'$. Let $m=x\i F(x)$. Then 
$m\in Z_H(g)$ and $m\in N_H(T)$. Let $\car_G,\car_H$ (resp. $\car'_G,\car'_H$)
be the set of roots of $Z_G^0(g),Z_H^0(g)$ (resp. $Z_G^0(g'),Z_H^0(g')$) with 
respect to $T$ (resp. $T'$). Then $F$ acts naturally on 
$\car_G,\car_H,\car'_G,\car'_H$. Let $\car'_G{}^+,\car'_H{}^+$ be a set of 
positive roots for $\car'_G,\car'_H$ respectively. It is well known that 

$\ep_{T'}\ep_{Z_G^0(g')}=
(-1)^{\sh(\a\in\car'_G{}^+;F(\a)\in\car'_G-\car'_G{}^+)}$,

$\ep_{T'}\ep_{Z_H^0(g')}=
(-1)^{\sh(\a\in\car'_H{}^+;F(\a)\in\car'_H-\car'_H{}^+)}$.
\nl
Now $\Ad(x)$ establishes bijections $\car_G\cong\car'_G$ and 
$\car_H\cong\car'_H$ under which the action of $F$ on $\car'_G,\car'_H$
corresponds to the action given by $\Ad(m)$ on $\car_G,\car_H$ (since $F$ acts
trivially on $\car_G,\car_H$). 
Also under this bijection, $\car'_G{}^+$ (resp. 
$\car'_H{}^+$) corresponds to a set of positive roots $\car_G^+$ (resp. 
$\car_H^+$) for $\car_G$ (resp. $\car_H$). It follows that

$\ep_{T'}\ep_{Z_G^0(g')}=
(-1)^{\sh(\a\in\car_G{}^+;\Ad(m)(\a)\in\car_G-\car_G{}^+)}$,

$\ep_{T'}\ep_{Z_H^0(g')}=
(-1)^{\sh(\a\in\car_H{}^+;\Ad(m)(\a)\in\car_H-\car_H{}^+)}$.
\nl
The right hand sides of the previous two equalities may be also interpreted as
the determinant of the linear map induced by $\Ad(m)$ on the group of 
characters of $T$. (We use that $m$ defines an element in the Weyl group of 
$Z_G^0(g)$ or $Z_H^0(g)$ with respect to $T$.) Hence those right hand sides 
coincide. It follows that (a) holds.

Now $\Ad(x)$ also induces an isomorphism $\bZ_H(g)\cong\bZ_H(g')$ under which 
the action of $F$ on $\bZ_H(g')$ corresponds to the action of $\Ad(m)$ on 
$\bZ_H(g)$ (since $F$ acts trivially on $\bar Z_H(g)$). Hence to prove (b) it 
is enough to show that $\Ad(m)$ acts trivially on $\bZ_H(g)$. Since 
$m\in Z_H(g)$, this follows from the well known fact that $\bZ_H(g)$ is 
commutative. The lemma is proved.

\subhead 3.18\endsubhead
In this subsection we assume that $(H,F)$ is like $(G,F)$ in 3.1. We assume
that $H=H_{der}$ and that $C$ is a semisimple conjugacy class in $H$. Let 
$a\in C$. We assume that $F(a)=a$ and that $F$ acts trivially on $\bZ_H(a)$. We
assume also that $F$ acts trivially on $Z=Z_{\tH}$.

(a) {\it Define $\mu:Z@>>>H^F/\pi(\tH^F)$ by $z\m\pi(h)$ where $h\in\tH$ 
satisfies $h\i F(h)=z$. Then $\mu$ is a group isomorphism.}
\nl
The proof is standard.

Let $Z'$ be the set of all $z\in Z$ such that $\ti a,z\ti a$ are conjugate in 
$\tH$; here $\ti a\in\pi\i(a)\sub\tH$. (This definition does not depend on the
choice of $\ti a$.) Note that $Z'$ is a subgroup of $Z$. Let $X$ be the set of
$H^F$-conjugacy classes contained in $C^F$.

(b) {\it Define $\mu':Z'@>>>X$ by $z\m\cc_{H^F}(\pi(g)a\pi(g)\i)$ where 
$g\in\tH$ is such that $g\i F(g)=g_1,g_1\ti ag_1\i=za$. Then $\mu'$ is a 
bijection.}
\nl
The proof is standard. 

(c) {\it The composition $Z'@>>>Z@>\mu>>H^F/\pi(\tH^F)$ 
(the first map is the inclusion) coincides with the
composition $Z'@>\mu'>>X@>>>H^F/\pi(\tH^F)$ (the second map is $A\m\pi(g)$ 
where $g\in A$).}
\nl
This follows from the definitions.

(d) {\it Assume that $N\ge 1$ is an integer such that $a^N=1$. Then
$z^N=1$ for all $z\in Z'$.}
\nl
Indeed if $z\in Z'$ we have $z\ti a=g\ti ag\i$ for some $g\in G$. Taking $N$-th
powers gives $z^N\ti a^N=g\ti a^Ng\i$ hence (using $\ti a^N\in Z$):
$z^N\ti a^N=gg\i\ti a^N=\ti a^N$, hence $z^N=1$.

Assume now that $(C'_2,C'_3,C'_5)$ is a $\fff$-triple for $H$. Let 
$a_n\in C'_n,n=2,3,5$. We assume that $F(a_n)=a_n$ and $F$ acts trivially on 
$\bZ_H(a_n)$ for $n=2,3,5$. 

Let $Z'_n,X_n$ be defined like $Z',X$ in terms of $C'_n,a_n$ instead of $C,a$.

(e) {\it The map $X_2\T X_3\T X_5@>>>H^F/\pi(\tH^F)$, given by
$(A,B,C)\m\pi(abc)$ where $(a,b,c)\in A\T B\T C$, is injective.}
\nl
Using (a),(b),(c), we see that this is equivalent to the statement that the map
$Z'_2\T Z'_3\T Z'_5@>>>Z$ given by multiplication in $Z$ is injective. From (d)
we see that any element $z\in Z'_n$ satisfies $z^n=1$. Then the injectivity of
the map above follows from the fact that an element of a finite abelian group 
can be written in at most one way as a product of elements of order dividing 
$2,3,5$.

In the remainder of this subsection we assume in addition that $H\in\cx'_{G,F}$
where $(G,F)$ as in 3.1 is such that $G=G_{ad}$, $G$ simple of type $E_m$,
$m=6,7,8$. (Recall that $H=H_{der}$.) 

(f) {\it The image $I$ of the map in (e) is exactly 
$U=\{x\in H^F/\pi(\tH^F);x^{30}=1\}$.}
\nl
By (d), we have $I\sub U$. Now $\x{I}=\x{Z'_2}\x{Z'_3}\x{Z'_5}$ and
$\x{U}$ can be easily determined by inspection of the various cases. We find 
that $\x{I}=\x{U}$ hence $I=U$.

((f) ought to be true also in type $D$, but we have not checked it, as we
don't need it.) In particular, the map in (e) is a bijection except if $(G,H)$
is of type $(E_8,D_5\T A_3)$ or $(E_7,A_3\T A_3\T A_1)$ in which case the image
has index $2$.

\proclaim{Lemma 3.19} Assume that $G=G_{ad}$. Let $F_0:G@>>>G$ be the Frobenius
map corresponding to an $\bF_p$-rational structure on $G$. 
There exists an integer $s_0\ge 1$ such that if $s\ge 1$ is an integer
divisible by $s_0$, and $H\in\cx'_{G,F_0^s}$ is semisimple, then  
$H$ is split over $\bF_{p^s}$.
\endproclaim
Let $G^*$ be as in the proof of 3.15. Let $F'_0:G^*@>>>G^*$ be the Frobenius 
map corresponding to an $\bF_p$-rational structure on $G^*$. As in the proof of
3.15, we see that it is enough to verify the following statement:

(a) There exists an integer $s_0\ge 1$ such that if $s\ge 1$ is an integer
divisible by $s_0$, and $g\in G^*{}^{F_0'{}^s}$ is a semisimple element 
such that $H=Z_{G^*}(g)$ is semisimple, then $H$ is split over $\bF_{p^s}$.
\nl
The set of all $g'\in G^*$ such that $g'$ is semisimple and $Z_{G^*}(g')$ is 
semisimple is a union of finitely many conjugacy classes $A_1,A_2,\do,A_m$ in
$G^*$. Pick $g_r\in A_r$ for $r\in[1,m]$. Let $H_r=Z_{G^*}(g_r)$. We can find 
an integer $s_0\ge 1$ such that $g_r\in G^*{}^{F_0'{}^{s_0}}$ and $H_r$ is 
split over $\bF_{p^{s_0}}$, $r=1,\do,m$. If $s\ge 1$ is an integer divisible by
$s_0$ then $g_r\in G^*{}^{F_0'{}^s}$ and $H_r$ is split over $\bF_{p^s}$, 
$r=1,\do,m$. Now let $g\in G^*{}^{F_0'{}^s}$ be such that $H=Z_{G^*}(g)$ is 
semisimple. We have $g\in A_r$ for some $r$. Hence $g\in A_r^{F_0'{}^s}$. Now
$Z_{G^*}(g_r)$ is connected since $G^*$ is simply connected. It follows that
$G^*{}^{F'{}^s}$ acts transitively (by conjugation) on $A_r^{F_0'{}^s}$. In 
particular, $g,g_r$ are conjugate under $G^*{}^{F'{}^s}$. It follows that
$H=Z_{G^*}(g)$ is split over $\bF_{p^s}$. The lemma is proved.

\head 4. Estimates\endhead
\proclaim{Lemma 4.1} Let $G,q,F$ be as in 3.1. Let $A,B,C$ be conjugacy classes
in $G^F$ such that $A\sub C_2,B\sub C_3,C\sub C_5$ where $(C_2,C_3,C_5)$ is a 
$\fff$-triple in $G$. Assume that the assumptions of 2.2 hold. Then 
$|\cg_{G,F;A,B,C}|\le\boc$ where $\boc$ is a constant depending only on $G$
(not on $F,q,A,B,C$). Moreover, if $(C_2,C_3,C_5)$ is not regular, then
$|\cg_{G,F;A,B,C}|\le\boc'q\i$ where $\boc'$ is a constant depending only on
$G$ (not on $F,q,A,B,C$).
\endproclaim
We use the notation of 2.1. We may assume that $(g_2,g_3,g_5)\in A\T B\T C$. 
Now $F$ acts naturally on $W$. Hence it acts naturally on the set of families 
of $W$. By \cite{\LU, 4.23} we have a partition 

$\ce_1(G^F)=\sqc_\cf\ce_{1,\cf}(G^F)$
\nl
($\cf$ runs over the $F$-stable families in $W$) such that the following holds:
if $\rh\in\ce_{1,\cf}(G^F)$, then 
$$f^\rh=\sum_Es_{\rh,E}f_E,$$
where $E$ runs through the representations in $\cf$ that are "$F$-stable",
$s_{\rh,E}$ are complex numbers that are bounded above (when $q$ varies) and
$f_E:\ct(G)^F@>>>\bc$ is as in 3.16. Using 3.16(a) we see that
$$R_{f^\rh,G}(1)=\sum_Es_{\rh,E}\sum_j\tr(\si_j,\Hom_W(E,\bv_j))q^j,$$
$$R_{f^\rh,Z_G(g_n)^0}(1)=\sum_Es_{\rh,E}\sum_j\tr(\si_j^{(n)},
\Hom_{W_n}(E|_{W_n},\bv_j^{(n)}))q^j,\quad n=2,3,5,$$
where $\si_j,\si_j^{(n)}$ are linear maps of finite order. Then

$R_{f^\rh,Z_G(g_2)^0}(1)R_{f^\rh,Z_G(g_3)^0}(1)R_{f^\rh,Z_G(g_5)^0}(1)$
\nl
is a linear combination with bounded coefficients of powers $q^j$ where 

$j\le b'_2(E)+b'_3(E')+b'_5(E'')$ with $E,E',E''\in\cf$
\nl
(hence $j\le\sum_{n=2,3,5}a'_n(\cf)$), 
while $R_{f^\rh,G}(1)$ is a linear combination 
with bounded coefficients of powers $q^j$ where $j\le b'(E)$ with $E\in\cf$ 
(hence $j\le a'(\cf)$) and in fact it is known that the coefficient of 
$q^{a'(\cf)}$ is non-zero and its inverse is bounded above. Now

$\cg_{G,F;A,B,C}=\sum_\cf\cg_{G,F;A,B,C,\cf}$
\nl
where $\cf$ runs over the $F$-stable families in $W$ and $\cg_{G,F;A,B,C,\cf}$ 
is defined like $\cg_{G,F;A,B,C}$ but with $\rh$ restricted to 
$\ce_{1,\cf}(G^F)$. The arguments above show, using 2.2(b), that each
$\cg_{G,F;A,B,C,\cf}$ is bounded above; and, in the case where $(C_2,C_3,C_5)$ 
is not regular, that each $q\cg_{G,F;A,B,C,\cf}$ is bounded above. The lemma 
follows.

\subhead 4.2\endsubhead
In the setup of 4.1, assume that $(C_2,C_3,C_5)$ is regular. By the method of 
4.1, we see that

(a) $q(\cg_{G,F;A,B,C}-1)$ is bounded above.

\proclaim{Lemma 4.3} Let $G$ be as in 3.1. Let $(C_2,C_3,C_5)$ be a regular 
$\fff$-triple in $G_{der}$. Then there exists $\ps\in\Hom(\ca_5,G)$ such that
$\ps$ gives rise to $(C_2,C_3,C_5)$ as in 1.7.
\endproclaim
Since this statement depends only on $G_{der}$ and
we can find an imbedding of $G_{der}$ into a connected reductive group with 
connected centre and derived subgroup $G_{der}$, we see that we
may assume that $G$ has connected centre. We may also assume that $\dim G>0$ 
and that the lemma is true
when $G$ is replaced by a group whose derived group has strictly smaller 
dimension than that of $G_{der}$. Assume that there is no $\ps$ as in the 
lemma (for $G$). Let $(g_2,g_3,g_5)\in C_2\T C_3\T C_5$. Since 
$g_n\in G_{der}$, we may choose $q,F$ as in 3.1 in such a way that 
$g_n\in\pi(\tG^F)$ for $n=2,3,5$. We may also assume that $G$,
$Z_G^0(g_2),Z_G^0(g_3),Z_G^0(g_5)$ are split over $\bF_q$ and that $q$ is 
large.

We write the identity in 3.13 for 
$(A,B,C)=(\cc_{G^F}(g_2),\cc_{G^F}(g_3),\cc_{G^F}(g_5))$. The left hand side of
that identity is $0$, by our assumption. We deduce that
$$\sum_{H\in\bar\cx'_{G,F}}\x{N_G(H)^F/H^F}\i X_{G,F;H;A,B,C}=0.\tag a$$
We show that

(b) $qX_{G,F;H;A,B,C}$ is bounded above (when $q,F$ vary) for any 
$H\in\bar\cx'_{G,F},H\ne G$.
\nl
With notation as in 3.12, we have  $qX_{G,F;H;A,B,C}=\sum'+\sum''$ where
$\sum',\sum''$ are given by
$$\sum_{(A',B',C')}t_{G,F;H;A,B,C;A',B',C'}\cg_{H,F;A',B',C'} 
qq^{-\dim Z^0_H}\sum_{\th\in H_\spa^F}\th(a'b'c');$$
in $\sum'$, (resp. $\sum''$), $(A',B',C')$ runs over all triples of conjugacy 
classes in $H^F$ such that $A'\sub A,B'\sub B,C'\sub C$ and
$A'\sub\ti A',B'\sub\ti B',C'\sub\ti C'$ where $(\ti A',\ti B',\ti C')$ is a 
$\fff$-triple in $H$, and 
$(\ti A',\ti B',\ti C')$ is regular (resp. non-regular) for $H$.

Now $\sum''$ is bounded above; indeed, in each term, $t_{G,F;H;A,B,C;A',B',C'}$
is bounded above (by 3.14), $q\cg_{H,F;A',B',C'}$ is bounded above (by 4.1) and
\linebreak $q^{-\dim Z^0_H}\sum_{\th\in H_\spa^F}\th(a'b'c')$ is bounded above
since $q^{-\dim Z^0_H}\x{H_\spa^F}$ is bounded above.

We show that $\sum'$ is bounded above. In each term, $t_{G,F;H;A,B,C;A',B',C'}$
is bounded above (by 3.14), $\cg_{H,F;A',B',C'}$ is bounded above (by 4.1)
and it is enough to show that 
$qq^{-\dim Z^0_H}\sum_{\th\in H_\spa^F}\th(a'b'c')$ is bounded above
or that
$$qq^{-\dim Z^0_H}\sum_{\th\in H_\he^F}\th(a'b'c')-
qq^{-\dim Z^0_H}\sum_{\th\in H_\he^F-H_\spa^F}\th(a'b'c')$$
is bounded above. Now 
$qq^{-\dim Z^0_H}\sum_{\th\in H_\he^F-H_\spa^F}\th(a'b'c')$ is bounded above 
since $qq^{-\dim Z^0_H}\x{H_\he^F-H_\spa^F}$ is bounded above (by 3.15(b)). 
Hence it is enough to show that $\sum_{\th\in H_\he^F}\th(a'b'c')=0$. Since 
$H^F_\he$ is a union of $H^F_*$-cosets in $H^F_\di$, it is enough to show that
$a'b'c'\n H_{der}$. Assume that $a'b'c'\in H_{der}$ for some (hence any)
$(a',b',c')\in A'\T B'\T C'$. Since $a',b',c'$ have orders dividing $2,3,5$, it
follows that each of $a',b',c'$ is in $H_{der}$. Hence $\ti A',\ti B',\ti C'$ 
are contained in $H_{der}$. Since $(\ti A',\ti B',\ti C')$ is regular in $H$, 
we may apply to it the induction hypothesis; we see that there exists
$\ps\in\Hom(\ca_5,H)$ such that $\ps(x_2)\in\ti A',\ps(x_3)\in\ti B'$,
$\ps(x_5)\in\ti C'$. Since $\ti A'\sub C_2,\ti B'\sub C_3,\ti C'\sub C_5$, we
see that $\ps(x_n)\in C_n$ for $n=2,3,5$. This contradicts our assumption. 
Thus, we have $a'b'c'\n H_{der}$ and the boundedness of $\sum'$ (hence (b)) is
established. 

Using now (a), we deduce that $qX_{G,F;G;A,B,C}$ is bounded above. Setting 
$t=t_{G,F;G;A,B,C;A,B,C},\cg=\cg_{G,F;A,B,C}$, we have 
$$qX_{G,F;G;A,B,C}=qt\cg q^{-\dim Z^0_G}\sum_{\th\in G_\di^F}\th(g_2g_3g_5)=
qt\cg q^{-\dim Z^0_G}\x{G_\di^F}$$
since $g_2g_3g_5\in\pi(\tG^F)$. From 3.14(b) we see that $t\i$ is bounded above
and it is clear that $q^{\dim Z^0_G}\x{G_\di^F}\i$ is bounded above. It follows
that $q\cg$ is bounded above. This contradicts 4.2. The lemma is proved.

\subhead 4.4\endsubhead
Assume now that $G$ (as in 1.2) is adjoint of type $E_8$ and that 
$(C_2,C_3,C_5)$ is the unique regular $\fff$ triple of $G$. Let $d$ be the
number of $G$-orbits on \linebreak
$\Hom_{reg}(\ca_5,G)$. To determine $d$ we may assume
that $k$ is as in 3.1. Let $q,F$ be as in 3.1. We write the identity 3.11(a) 
for $(A,B,C)=(C_2^F,C_3^F,C_5^F)$ (these are three conjugacy classes of $G^F$).
Assume that we can evaluate the right hand side of 3.11(a) for infinitely many
$q$ and that it is of the form $\De q^D+$ lower powers of $q$ where $\De$ is a
constant. Then 3.11(a) implies that $\{(a,b,c)\in C_2\T C_3\T C_5;abc=1\}$ has
exactly $\De$ irreducible components of dimension $D$ and it follows that 
$d=\De$. (This strategy was suggested in \cite{\SE}.)

\proclaim{Theorem 4.5} Let $G$ be as in 1.2 and let $(C_2,C_3,C_5)$ be a
regular $\fff$-triple of $G$. Assume that $G=G_{ad}$. Let $N=\x{Z_{\tG}}$. Let
$g_n\in C_n$ and let $N_n=\x{\bZ_G(g_n)}$, $n=2,3,5$. Then, up to 
$G$-conjugacy, there are exactly $\fra{N}{N_2N_3N_5}$ elements 
$\ps\in\Hom_{reg}(\ca_5,G)$ that give rise to $(C_2,C_3,C_5)$.
\endproclaim
We may assume that $G$ is simple. In the case where $G$ is of type $A$ or $D$,
the theorem can be obtained from the results in \S1. In the rest of the proof
we assume that $G$ is of type $E_m,m=6,7,8$ (although the same proof should
work without this assumption).

We will carry out the strategy in 4.4, using the results in \S2, \S3. We may 
assume that $k$ is as in 3.1. Let $F_0$ be a Frobenius map for an 
$\bF_p$-rational structure on $G$. Let $q=p^{s'},F=F_0^{s'}$ where $s'$ is 
sufficiently large and divisible by a fixed integer $s'_0\ge 1$.
By choosing $s'_0$ apropriately we may assume that $F$ acts trivially on 
$Z_{\tG}$ that $g_n\in\pi(\tG^F)$, that
$Z_G^0(g_n)$ is split over $F_q$ and that $F$ acts trivially on 
$\bar Z_G(g_n)$ for $n=2,3,5$. (The choice of $s'_0$ will be made more 
precise later in the proof.)
Let 
$\de=3\nu-\nu_2-\nu_3-\nu_5$. (Notation of 2.1.) Let $\De=\fra{N}{N_2N_3N_5}$. 
We will show that
$$q(q^{-\de}\sh\{(g_2,g_3,g_5)\in C^F_2\T C^F_3\T C^F_5;g_2g_3g_5=1\}-\De)
\tag a$$
is bounded above when $q$ varies or equivalently, that
$$q(q^{-\de}\sum_{(A,B,C)\in\cx}
\sh\{(g_2,g_3,g_5)\in A\T B\T C;g_2g_3g_5=1\}-\De)$$
is bounded above. Here $\cx$ is the set of triples $(A,B,C)$ of conjugacy 
classes in $G^F$ such that $A\sub C_2,B\sub C_3,C\sub C_5$. By 3.13, it is 
enough to show that

(b) $q(\sum_{(A,B,C)\in\cx}X_{G,F;G;A,B,C}-\De)$ is bounded above,
\nl 
(c) $q\sum_{(A,B,C)\in\cx}X_{G,F;H;A,B,C}$ is bounded above for any 
$H\in\bar\cx'_{G,F},H\ne G$.
\nl
We have

$\sum_{(A,B,C)\in\cx}X_{G,F;G;A,B,C}=
\sum_{A,B,C}t_{G,F;G;A,B,C;A,B,C}\cg_{G,F;A,B,C}\sum_{\th\in G_\di^F}\th(abc)$
\nl
where $(a,b,c)\in A\T B\T C$. The sum over $\th$ is $0$ unless 
$abc\in\pi(\tG^F)$ (or equivalently, each of $a,b,c$ is in $\pi(\tG^F)$), in 
which case it is $N$. Thus, (b) is equivalent to the statement that

$q(t_{G,F;G;A,B,C;A,B,C}\cg_{G,F;A,B,C}N-\De)$  is bounded above;
\nl
here $(A,B,C)\in\cx$ is uniquely determined by the condition

$A\sub\pi(\tG^F),B\sub\pi(\tG^F),C\sub\pi(\tG^F)$.
\nl
This follows from the fact that $q(t_{G,F;G;A,B,C;A,B,C}-(N_2N_3N_5)\i)$ is 
bounded above (see 3.14, 3.17) and 
$q(\cg_{G,F;A,B,C}-1)$ is bounded above (see 4.2). Thus, (b) holds.

Now, in (c), we have (as in the proof of 4.3)

$q\sum_{(A,B,C)\in\cx}X_{G,F;H;A,B,C}=\sum'+\sum''$ 
\nl
where $\sum',\sum''$ are given by 
$$\sum_{(A',B',C')}t_{G,F;H;A,B,C;A',B',C'}\cg_{H,F;A',B',C'}
qq^{-\dim Z^0_H}\sum_{\th\in H_\spa^F}\th(a'b'c');$$
in $\sum'$, (resp. $\sum''$), $(A',B',C')$ runs over the set $\cx'$ (resp.
$\cx''$) consisting of all triples of conjugacy 
classes in $H^F$ such that $A'\sub\ti A',B'\sub\ti B',C'\sub\ti C'$ where 
$(\ti A',\ti B',\ti C')$ is a $\fff$-triple in $H$, 
$\ti A'\sub C_2,\ti B'\sub C_3,\ti C'\sub C_5$ and $(\ti A',\ti B',\ti C')$ is
regular (resp. non-regular) for $H$; moreover $(A,B,C)$ is the unique triple of
conjugacy classes in $G^F$ such that $A'\sub A,B'\sub B,C'\sub C$.

Now $\sum''$ is bounded above (exactly as in the proof of 4.3).

We now estimate $\sum'$.
 
Assume first that $\dim Z_H>0$. Again, $\sum'$ is bounded above: in each term,
$t_{G,F;H;A,B,C;A',B',C'}$ is bounded above (by 3.14), $\cg_{H,F;A',B',C'}$ is 
bounded above (by 4.1)  and it is enough to show that 
$qq^{-\dim Z^0_H}\sum_{\th\in H_\spa^F}\th(a'b'c')$ is bounded above or that

$qq^{-\dim Z^0_H}\sum_{\th\in H_\he^F}\th(a'b'c')-
qq^{-\dim Z^0_H}\sum_{\th\in H_\he^F-H_\spa^F}\th(a'b'c')$
\nl
is bounded above. Now
$qq^{-\dim Z^0_H}\sum_{\th\in H_\he^F-H_\spa^F}\th(a'b'c')$ is bounded above 
since $qq^{-\dim Z^0_H}\x{H_\he^F-H_\spa^F}$ is bounded above (by 3.15(b)). 
Hence it is enough to show that $\sum_{\th\in H_\he^F}\th(a'b'c')=0$. Since 
$H^F_\he$ is a union of $H^F_*$-cosets in $H^F_\di$, it is enough to show that
$a'b'c'\n H_{der}$. Assume that $a'b'c'\in H_{der}$ for some (hence any)
$(a',b',c')\in A'\T B'\T C'$. Since $a',b',c'$ have orders dividing $2,3,5$, it
follows that each of $a',b',c'$ is in $H_{der}$. Hence $\ti A',\ti B',\ti C'$
are contained in $H_{der}$. Since $(\ti A',\ti B',\ti C')$ is regular in $H$, 
we may apply to it Lemma 4.3; we see that there exists $\ps\in\Hom(\ca_5,H)$ 
such that $\ps(x_2)\in\ti A'$, $\ps(x_3)\in\ti B'$, $\ps(x_5)\in\ti C'$. Since
$\ti A'\sub C_2,\ti B'\sub C_3,\ti C'\sub C_3$, we see that $\ps(x_n)\in C_n$ 
for $n=2,3,5$. Since $(C_2,C_3,C_5)$ is regular for $G$, we have 
$\ps\in\Hom_{reg}(\ca_5,G)$. This contradicts the fact that $Z_G(\ps(\ca_5))$ 
contains the non-trivial torus $Z_H^0$. Thus, we have $a'b'c'\n H_{der}$ and 
the boundedness of $\sum'$ is established. 

Next, assume that $\dim Z_H=0$. We show that $\sum'$ is bounded above. Since 
for $(A',B',C')\in\cx'$,
$q(\cg_{H,F;A',B',C'}-1)$ is bounded above (see 4.2) 
and \linebreak
$t_{G,F;H;A,B,C;A',B',C'}$ is bounded above (see 3.14), it is enough to 
show that, for any regular $\fff$-triple $(C'_2,C'_3,C'_5)$ of $H$ such that
$C'_n\sub C_n$, $n=2,3,5$, and for any $\th\in H_\spa^F$,

$\sum_{(A',B',C')\in\cy}qt_{G,F;H;A,B,C;A',B',C'}\th(a'b'c')$
\nl
is bounded above. Here 
$\cy=\{(A',B',C')\in\cx'; A'\sub C'_2,B'\sub C'_3,C'\sub C'_5\}$ and
$(a',b',c')\in A'\T B'\T C'$. Using 3.14(c) we see that it is enough to show 
that, for any $\th\in H_\spa^F$,
 
$\sum_{(A',B',C')\in\cy}\t_{A',B',C'}\th(a'b'c')=0$
\nl
where
$$\t_{A',B',C'}=\fra{\ep_{Z^0_G(a')}\ep_{Z^0_H(a')}
\ep_{Z^0_G(b')}\ep_{Z^0_H(b')}\ep_{Z^0_G(c')}\ep_{Z^0_H(c')}}
{\x{\bZ_H(a')^F}\x{\bZ_H(b')^F}\x{\bZ_H(c')^F}}.$$
By choosing $s'_0$ appropriately, we may assume that $H$ is split over 
$\bF_q$ (see 3.19), that $F$ acts trivially on $Z_{\tH}$, that $C'_n$ 
contains an 
$\bF_q$-rational point $h_n$ such that $Z_H^0(h_n)$ is split over $F_q$ and $F$
acts trivially on $\bZ_H(h_n)$, $n=2,3,5$. Using 3.17, we see then that, in the
last fraction, the numerator is $1$ and the denominator is
${\x{\bZ_H(h_2)}\x{\bZ_H(h_3)}\x{\bZ_H(h_5)}}$; in particular,
$\t_{A',B',C'}$ is independent of $A',B',C'$. Thus, we are reduced to showing
that

$\sum_{(A',B',C')\in\cy}\th(a'b'c')=0$
\nl
Using the results at the end of 3.18, we see that the last equality is
equivalent to the equality

$\sum_{x\in H^F/\pi(\tH^F);x^{30}=1}\th(x)=0$
\nl
where $\th$ is regarded as a character of $H^F/\pi(\tH^F)$. More precisely, if
$(G,H)$ is not of type $(E_8,D_5\T A_3)$ or $(E_7,A_3\T A_3\T A_1)$, we must
show that $\sum_{x\in H^F/\pi(\tH^F)}\th(x)=0$ and this follows from the fact 
that $\th\ne 1$ (since $\th\in H^F_\spa$). If $(G,H)$ is of type 
$(E_8,D_5\T A_3)$ or $(E_7,A_3\T A_3\T A_1)$, we must show that the restriction
of $\th$ to the subgroup of $H^F/\pi(\tH^F)$ ($\cong\bz/4\bz$ or 
$\cong\bz/4\bz\opl\bz/2\bz$) consisting of all $x$ such that $x^2=1$ is 
$\ne 1$; this again follows from $\th\in H^F_\spa$.

This proves (c). Thus, (a) is proved.

Let $\bx=\{(a,b,c)\in C_2\T C_3\T C_5;abc=1\}$ (an algebraic variety defined 
over $\bF_q$). From (a) we see that there exists an integer $M>0$ such that  
$q^s(q^{-\de s}\x{\bx^{F^s}}-\De)\le M$ for $s=1,2,\do$. As it is known, this
implies that among the irreducible components of $\bx$ there are exactly $\De$
of maximum dimension, and that maximum dimension is $\de$. By 1.3, any 
$G$-orbit on $\bx$ has dimension $\de$. It follows that $\bx$ has exactly $\De$
$G$-orbits. Since $\bx$ is naturally in bijection with the set of all 
$\ps\in\Hom(\ca_5,G)$ that give rise to $(C_2,C_3,C_5)$ we see that the
theorem is proved.

\subhead 4.6\endsubhead
For $G=G_{ad}$ of type $E_6,E_7,E_8$, the fraction $\fra{N}{N_2N_3N_5}$ is
$\fra{3}{1\cdot 3\cdot 1}=1$, $\fra{2}{2\cdot 1\cdot 1}=1$,
$\fra{1}{1\cdot 1\cdot 1}=1$ respectively. We see that {\it in type $E_8$
there is exactly one $\ps\in\Hom_{reg}(\ca_5,G)$ up to $G$-conjugacy}.

For $G=G_{ad}$ of type $A_1,A_2,A_3,A_4,A_5$, the fraction $\fra{N}{N_2N_3N_5}$
is $\fra{2}{2\cdot 1\cdot 1}=1$, $\fra{3}{1\cdot 3\cdot 1}=1$,
$\fra{4}{2\cdot 1\cdot 1}=2$, $\fra{5}{1\cdot 1\cdot 5}=1$,
$\fra{6}{2\cdot 3\cdot 1}=1$ respectively. 

For $G=G_{ad}$ of type $D_4,D_5,D_6,D_8$, the fraction $\fra{N}{N_2N_3N_5}$ is
$\fra{4}{4\cdot 1\cdot 1}=1$, $\fra{4}{2\cdot 1\cdot 1}=2$,
$\fra{4}{4\cdot 1\cdot 1}=1$, $\fra{4}{4\cdot 1\cdot 1}=1$ respectively. 

\head 5. Relation to homomorphisms $PGL_2(k)@>>>G$\endhead
\subhead 5.1\endsubhead
In this section we assume that either $p=0$ or $p$ is large enough. Let $G$ be
as in 1.2. Let $\Phi:PGL_2(k)@>>>G$ be a homomorphism. Composing $\Phi$ with 
one of the two regular homomorphisms $\ca_5@>>>PGL_2(k)$ (up to 
$PGL_2(k)$-conjugacy, see 1.8), we obtain two homomorphisms $\ca_5@>>>G$ 
(which may or may not be $G$-conjugate).

In \cite{\SE}, it is pointed out that by applying this construction to a $\Phi$
whose image contains a regular unipotent element in $G$ (of type $E_m$) we 
obtain a regular $\ps$. Applying the same procedure for $G=G_{ad}$ of type 
$A_{m-1}$, we see that all regular $\ps$ in 1.8 can be thus obtained for
$m=2,3,5,6$ but only one of the two regular $\ps$ is thus obtained for $m=4$.
This procedure can be applied for $G$ as in 1.9. In that case, for $m=4$, the 
first two $\ps$ in 1.9(a) are obtained; for $m=5$, the first $\ps$ in 1.9(b) is
obtained; for $m=6$, the first $\ps$ in 1.9(c) is obtained; for $m=8$, the
first $\ps$ in 1.9(d) is obtained.

Let us analyze what happens if we use some other $\Phi$.

Let $G$ be as in 1.9 and assume that $\Phi$ is such that its image contains a
subregular unipotent element of $G$. For $m=4$, we thus obtain the third and 
fourth $\ps$ in 1.9(a); for $m=5$, we thus obtain the second $\ps$ in 1.9(b);
for $m=6$, we thus obtain the second and third $\ps$ in 1.9(c).

Let $G$ be as in 1.9(d). If $\Phi$ is such that its image contains a unipotent
element whose Jordan blocks in the standard representation have sizes
$1,3,5,7$, then the associated $\ps$ is regular.

If $G$ (as in 1.10) is of type $E_8$ and $\Phi$ is such that the centralizer of
its image is the symmetric group in $5$ letters, then the $\ps$ obtained from
$\Phi$ by the procedure above is again regular (the corresponding $\fff$-triple
is the one described in 1.10). Using 4.5, the $\ps$ thus obtained must be the
same up to conjugacy as the one attached to a regular unipotent class.

These arguments together with 4.5 show that, if $G=G_{ad}$ is simple, then any
$\ps\in\Hom_{reg}(\ca_5,G)$ is obtained from some $\Phi:PGL_2(k)@>>>G$, except
for one case in type $A_3$.

\enddocument